\input miniltx   
\def\Gin@driver{pdftex.def}
\input graphicx.sty
\resetatcatcode

 at 21truept
\font\eightrm=cmr8
\font\twelverm=cmr12
\font\tenrm=cmr10

\font\eighteenbf=cmbx14 scaled\magstep1
\font\twelvebf=cmbx10 scaled\magstep1

\topskip .5in
\leftskip 2 cm
\hsize=5.75 in
\vsize=8.5in

\settabs 1 \columns
\twelverm

\centerline{ \eighteenbf  Inverse Quadratic Transportation Problem}
\bigskip
\bigskip
\tenrm
\settabs 1 \columns
\+\hfill { Afrooz Jalilzadeh}\hfill&\cr 
\medskip
\eightrm
\+\hfill Tehran University, Mathematics Department \hfill&\cr
\medskip
\tenrm
\+\hfill Erfan Yazdandoost Hamedani\hfill&\cr
\medskip
\eightrm
\+\hfill Tehran University, Mathematics Department \hfill&\cr
\bigskip
\bigskip

\hsize=5.45in
\leftskip 2cm
\noindent {\bf ABSTRACT.} Many research has been conducted about quadratic programming and inverse optimization. In this paper we present the combination aspect of these subjects, applying on transportation problem. First, we obtain the inverse form of quadratic tranportation problem under $L_1$ norm by using duality as well as introducing the optimal value. Then, we do the same process for inverse quadratic transportation problem (IQTP) under $L_\infty$ norm.

\hsize=5.75in
\leftskip 1cm
\tenrm
\bigskip
\bigskip

\noindent {\twelvebf 1. Introduction.}
\bigskip
\noindent The classical transportation problem (TP) is a well- structured problem which has been studied by many researchers. Consider $m$ origin points, where origin $i$ has a supply of $s_i$ units of a particular item (commodity). Additionally, suppose that there are $n$ destinations points, where destination $j$ requires $d_j$ units of the commodity. We assume that $s_i$  ,$d_j\ge0$ for all $i,j$. There is a unit cost $c_{ij}$ for each link $(i,j)$, from origin $i$ to destination $j$. The problem is to determine a feasible “shipping pattern” from origins to destinations that minimizes the total transportation cost. Let $x_{ij}$ be the number of units shipped along link $(i, j)$ from origin $i$ to destination $j$ [B1-chapter 10]. Without loss of generality we can assume that the TP is balanced which means the total supplies equals the total demands. Therefore, The balanced transportation problem can be represented as follows:
$${\rm Minimize} \ \ \sum_{i=1}^n \sum_{j=1}^m c_{ij}x_{ij} $$
$$s.t. \ \ \ \ \ \sum_{j=1}^m x_{ij}= s_i, \ \ \ \ \ \ \  i=1,...,n $$
$$\ \ \ \ \ \ \ \ \ \sum_{i=1}^n x_{ij}= d_j, \ \ \ \ \  j=1,...,m $$
$$x_{ij}\ge 0 , \ \ \ \ \  i=1,...,n, \ \ \ \  j=1,...,m $$

\par Beside the TP, quadratic programs are an important class of problems on their own and as subproblems in methods for general optimization problems. The quadratic programming problem involves minimizing a quadratic function subject to linear constraints. A general formulation is

$$ {\rm Minimize}\ {1 \over 2} x^T Q x + c^T x$$
$$s.t. \ a_i^T x=b_i, \ \ i=1,...,N  $$
$$\ \ \ \ a_i^T x \ge b_i, \ \ i=n+1,...,N+M $$
\noindent where $Q \in R^{(N+M) \times (N+M)}$ is symmetric. When $Q$ is semipositive definite matrix, the objective function is convex; there are polynomial-time algorithms [K1] to solve these relatively easy problems. Furthermore, Bomze and Danninger [B2] introduced a finite algorithm for solving general quadratic problems, which consist of maximizing a non-concave quadratic function over a polyhedron in n-dimensional Euclidean space. In addition, there are some research about integer and 0-1 quadratic programming by Bretthauer, Shetty and Syam [B3, B4] and by Chaovalitowangse, Pardalas and Prokopyev [C1] respectively. There are limited research which have presented algorithms for solving transportation problems with quadratic function cost coefficients. Adlakha and Kowalski [A1] provided the absolute point algorithm for solving quadratic transportation problem (QTP), Cosares and Hochbaum [C2] presented strongly polynomial algorithms for a QTP with a fixed number of sources.

\par Another most recent and important problem in operation research is inverse optimization. In the past few decades, there has been an interest in inverse optimization problem in the operation research community, and many aspect of these problems have been studied. The inverse problem comprises change of cost vector such that a feasible solution of the primal problem becomes optimal with respect to minimum change in the cost vector. Zang and Liu [Z2] studied inverse assignment and minimum cost flowproblems; Yang et al. [Y1] and Zhang and Cai [Z1] have studied the inverse minimum cut problems; Xu and Zhang [X1] have studied the inverse shortest path problem; and Ahuja and Orlin [A2] provided various references in the area of inverse optimization. Some research have been studied on inverse of quadratic programming as well. Zhang J., Zhang Li. and Xiao [Z3] have studied a perturbation approach for an inverse quadratic programming problem and recently, Jain and Arya [J1] presented inverse optimization for quadratic programming problems. 

\par Our purpose in this paper is to obtain the inverse form of QTP under $L_1$ and $L_\infty$ norm. In particular, by using some substitution the formulation of the problem under these norms becomes linear. In addition, we introduce the optimal value of the inverse quadratic transportation problem (IQTP).

\bigskip
\par \noindent {\twelvebf 2. Formulating Inverse of Quadratic Transportation problem.} 
\bigskip
\noindent In the real world the cost function of a TP can be represented by a quadratic form, so it seems that more research about QTP in different aspects could be useful. In QTP we have $n$ suppliers which can ship units of a product to any of the $m$ demands where the objective function is quadratic. We can suppose the supplies and demands are non-negative and the problem is balanced. The problem is to minimize the total transportation cost within limited supplies, the demands are satisfied. QTP is formulated as follows:

$${\rm Minimize}\ {1 \over 2} x^T Q x + c^T x \eqno(1a)$$
$$s.t. \ \ \ \ \ \sum_{j=1}^m x_{ij}=s_i, \ \ \ \ \ \  i=1,...,n \eqno(1b)$$
$$\ \ \ \ \ \ \ \ \ \sum_{i=1}^n x_{ij}= d_j, \ \ \ \ \  j=1,...,m \eqno(1c)$$
$$\ \ \ \ \ \ \ \ \ \ x_{ij}\ge 0 \eqno(1d)$$

\par Now we try to study the inverse form of quadratic transportation problem (IQTP) under $L_1$ and $L_{\infty}$ norm. The inverse problem comprises change of cost vector such that a feasible solution of the primal problem becomes optimal with respect to minimum change in the cost vector. Consider the following general convex QTP:
$${\rm Minimize} \ {1\over 2}x^T Qx+c^Tx \eqno(2a)$$
$$s.t. \ \ \ Ax= b \eqno(2b)$$
$$\ \ \ \ \ \ x\ge 0 \eqno(2c)$$
where $Q$ is a symmetric matrix, elements of matrix $A$ are cofficients of $x_{ij}$ and $b$ is the right hand side vector in constraints of (1).
\par In the objective function of this problem, the cost consists of the matrix $Q$ and the vector $c$. Consider the (2) with new cost, symmetric matrix $H$ and vector $d$:
$${\rm Minimize} \ {1\over 2}x^T Hx+d^Tx \eqno(3a)$$
$$s.t. \ \ \ Ax= b \eqno(3b)$$
$$\ \ \ \ \ \ x\ge 0 \eqno(3c)$$

\par For obtaining the inverse of QTP, we shall consider the mathematical formula of general inverse problem for the QTP. Suppose $x^0$ is a given feasible solution of the problem (2). We intend to make $x^0$ be the optimal solution for the problem (3) with minimum change. The objective function of IQTP is sum of subtraction of matrices cost and vectors cost. Thus, it is sum of a matrix norm and a vector norm. Hence, the IQTP is as follows:
$${\rm Minimize}\ \| Q-H \| + \| c-d \| \eqno(4a)$$
$$x^0\in O_d \eqno(4b)$$
$$H\ is\ symmetric \eqno(4c)$$
where $O_d$ is the set of optimal solutions of the (3). The condition (4b) can be restated using duality of the problem (3). Therefore, according to the Kuhn Tucher and complementary slackness conditions, the condition (4b) is equivalent to:
$$x^{0T}H+d^T-w^1A+w^2=0 \eqno(5a)$$
$$w^2x^0=0 \eqno(5b)$$
$$w^2\ge 0 \eqno(5c)$$
where $w^1 \in R^{m+n}$ and $w^2\in R$ are lagrangian multipliers.
\par For preventing from confusion in indexes we shall consider the elements of the matrices, such as $Q_{(mn)\times (mn)}$, in the form of $Q_{(kl)(ij)}$ where $k,i\in \{1,...,n\}$ and $l,j\in \{1,...,m\}$. We consider each row and each column of $Q$ as $m$ parts where each part has $n$ elements. Thus, $(kl)$ in $Q_{(kl)(ij)}$ shows that we consider the $l$th element of $k$th part in rows. Similarly, $(ij)$ in $Q_{(kl)(ij)}$ shows that we consider the $j$th element of $i$th part in columns. This notation helps to write the multiplication of these matrices to vectors in a simpler way. Let $J=\{(i,j)| 1\le i\le n,\ 1\le j \le m\}$, $L=\{(i,j)\in J | \ x^0_{ij}=0 \}$ and $F=\{(i,j)\in J | \ x^0_{ij}>0 \}$. By applying the conditions (5a)-(5c) to the problem (4) we obtain the IQTP:
\bigskip
$${\rm Minimize}\ \| Q-H \| + \| c-d \| \eqno(6a)$$
$$s.t.\ w^1_i+w^1_{n+j}-\sum_{l=1}^m\sum_{k=1}^n H_{(kl)(ij)}x^0_{kl}=d_{ij},\ (i,j)\in F \eqno(6b)$$
$$\ w^1_i+w^1_{n+j}-\sum_{l=1}^m\sum_{k=1}^n H_{(kl)(ij)}x^0_{kl}+w^2_{ij}=d_{ij},\ (i,j)\in L \eqno(6c)$$
$$H_{(kl)(ij)}=H_{(ij)(kl)}, \ (i,j),(k,l)\in J \eqno(6d)$$
$$w^2_{ij}\ge 0,\ (i,j)\in J \eqno(6e)$$

\bigskip
\par \noindent {\twelvebf 3. Inverse Quadratic Transportation Problem Under $L_1$ norm.}
\bigskip

\noindent In this section, we will consider the IQTP problem under $L_1$ norm. Therefore, our purpose is to study the following problem:

$${\rm Minimize}\ \| Q-H \|_1 + \| c-d \|_1 \eqno(7a)$$
$$s.t.\ w^1_i+w^1_{n+j}-\sum_{l=1}^m\sum_{k=1}^n H_{(kl)(ij)}x^0_{kl}=d_{ij},\ (i,j)\in F \eqno(7b)$$
$$\ w^1_i+w^1_{n+j}-\sum_{l=1}^m\sum_{k=1}^n H_{(kl)(ij)}x^0_{kl}+w^2_{ij}=d_{ij},\ (i,j)\in L \eqno(7c)$$
$$H_{(kl)(ij)}=H_{(ij)(kl)}, \ (i,j),(k,l)\in J \eqno(7d)$$
$$w^2_{ij}\ge 0,\ (i,j)\in J \eqno(7e)$$

\par $L_1$ norm of a matrix is the maximum absolute column sum of the matrix, so $\| H-Q \|_1={\rm max}_{(i,j)\in J}(\sum_{(k,l)\in J} |H_{(kl)(ij)}-Q_{(kl)(ij)}|)$ which we denote it by $\theta$. Note that that the objective function of (6) is nonlinear, so by using substitution of $d_{ij}-c_{ij}=\alpha_{ij}-\beta_{ij},\ |d_{ij}-c_{ij}|=\alpha_{ij}+\beta_{ij}$ and $H_{(kl)(ij)}-Q_{(kl)(ij)}=\Gamma_{(kl)(ij)}-\Delta_{(kl)(ij)},\ |H_{(kl)(ij)}-Q_{(kl)(ij)}|=\Gamma_{(kl)(ij)}+\Delta_{(kl)(ij)}$ where $\alpha_{ij},\ \beta_{ij},\ \Gamma_{(kl)(ij)},\ \Delta_{(kl)(ij)}\ge 0$, $\Gamma$ and $\Delta$ are symmetric matrices, we obtain the following problem which is linear form of IQTP and equivalent to the problem (7).
$${\rm Minimize}\ \theta+\sum_{i=1}^n\sum_{j=1}^m (\alpha_{ij}+\beta_{ij}) \eqno(8a)$$
$$s.t.\ \sum_{l=1}^m\sum_{k=1}^n [\Delta_{(kl)(ij)}-\Gamma_{(kl)(ij)}]x^0_{kl}-\alpha_{ij}+\beta_{ij}=c^{\pi}_{ij},\ (i,j)\in F \eqno(8b)$$
$$\ \sum_{l=1}^m\sum_{k=1}^n [\Delta_{(kl)(ij)}-\Gamma_{(kl)(ij)}]x^0_{kl}-\alpha_{ij}+\beta_{ij}+w^2_{ij}=c^{\pi}_{ij},\ (i,j)\in L \eqno(8c)$$
$$\Gamma_{(kl)(ij)}=\Gamma_{(ij)(kl)}, \ (i,j),(k,l)\in J \eqno(8d)$$
$$\Delta_{(kl)(ij)}=\Delta_{(ij)(kl)}, \ (i,j),(k,l)\in J \eqno(8e)$$
$$\theta \ge \sum_{l=1}^m\sum_{k=1}^n [\Delta_{(kl)(ij)}+\Gamma_{(kl)(ij)}],\ (i,j)\in F\cup L \eqno(8f)$$
$$\Gamma_{(kl)(ij)},\ \Delta_{(kl)(ij)},\ \alpha_{ij},\ \beta_{ij},\ w^2_{ij}\ge 0,\ (i,j),(k,l)\in J \eqno(8g)$$
\smallskip
where $c^{\pi}_{ij}=c_{ij}-(w^1_i+w^1_{n+j})+\sum_{l=1}^m\sum_{k=1}^n (Q_{(kl)(ij)}x^0_{kl})$.
\bigskip
We will now simplify equations (8b) and (8c). At first we note that in an optimal solution of (8), $\alpha_{ij}$ and $\beta_{ij}$ cannot take positive values and similarly $\Delta_{(kl)(ij)}$ and $\Gamma_{(kl)(ij)}$ cannot take positive values because if we can reduce both by small amount $\epsilon >0$ the constraints are still satisfied while the objective function value is strictly improved. Now by placing some conditions on $c_{ij}^\pi$ and using the aforementioned discussion we can obtain the optimal value of IQTP. There are two cases.
\medskip

\noindent {\bf Case 1.} $(i,j)\in F$. (i) If $c_{ij}^{\pi}>0$ and $\beta_{ij}>0$, then $\alpha_{ij}=\Gamma_{(kl)(ij)}=\Delta_{(kl)(ij)}=0$, so $c_{ij}^{\pi}=\beta_{ij}$. (ii)  If $c_{ij}^{\pi}>0$ and $\beta_{ij}=0$, then $\alpha_{ij}=\Gamma_{(kl)(ij)}=0$ and there exist an index $(kl)$ such that $x_{kl}^0>0$, because otherwise we obtain $c_{ij}^{\pi}=0$ which is contradiction to the assumption. Therefore, there exist an index $(kl)$ such that $\Delta_{(kl)(ij)}x_{kl}^0=c_{ij}^{\pi}$, so $\Delta_{(kl)(ij)}={|c_{ij}^{\pi}| \over x_{kl}^0}$. (iii) If $c_{ij}^{\pi}<0$ and $\alpha_{ij}>0$, then $\beta_{ij}=\Gamma_{(kl)(ij)}=\Delta_{(kl)(ij)}=0$, so $\alpha_{ij}=-|c_{ij}^{\pi}|$. (iv)  If $c_{ij}^{\pi}<0$ and $\alpha_{ij}=0$, then $\beta_{ij}=\Delta_{(kl)(ij)}=0$ and similar to the (ii) we have an index $(kl)$ such that $-\Gamma_{(kl)(ij)}x_{kl}^0=c_{ij}^{\pi}$, so $\Gamma_{(kl)(ij)}={|c_{ij}^{\pi}| \over x_{kl}^0}$. (v)  If $c_{ij}^{\pi}=0$, then $\alpha_{ij}=\beta_{ij}=\Gamma_{(kl)(ij)}=\Delta_{(kl)(ij)}=0$.
\medskip
\noindent {\bf Case 2.} $(i,j)\in L$. (i) If $c_{ij}^{\pi}>0$, then $\alpha_{ij}=\beta_{ij}=\Gamma_{(kl)(ij)}=\Delta_{(kl)(ij)}=0$ and $c_{ij}^{\pi}=w_{ij}^2$. (ii) If $c_{ij}^{\pi}<0$ and $\alpha_{ij}>0$, then $w_{ij}^2=\beta_{ij}=\Gamma_{(kl)(ij)}=\Delta_{(kl)(ij)}=0$, so $\alpha_{ij}=-|c_{ij}^{\pi}|$. (iii) If $c_{ij}^{\pi}<0$ and $\alpha_{ij}=0$, then $w_{ij}^2=\beta_{ij}=\Delta_{(kl)(ij)}=0$, so again similar to the (iv) of case 1, we have $\Gamma_{(kl)(ij)}={|c_{ij}^{\pi}| \over x_{kl}^0}$. (iv) If $c_{ij}^{\pi}=0$, then $w_{ij}^2=\alpha_{ij}=\beta_{ij}=\Gamma_{(kl)(ij)}=\Delta_{(kl)(ij)}=0$.
\medskip
\par According to the cases the optimal values of IQTP are: 
$$d_{ij}^*=\cases {c_{ij}+|c_{ij}^\pi| \ &if $\ c_{ij}^\pi < 0, \ (i,j)\in F\cup L , \ \alpha_{ij}>0$ \cr \cr  c_{ij}-|c_{ij}^\pi| \ &if $\ c_{ij}^\pi > 0, \ (i,j)\in F, \ \beta_{ij}>0$ \cr \cr c_{ij} \ & \ otherwise \cr}$$
\smallskip

$$H_{(kl)(ij)}^*=\cases {Q_{(kl)(ij)}+{|c_{ij}^\pi| \over x_{kl}^0} \ &if $\ c_{ij}^\pi < 0, \ (i,j)\in F\cup L , \ \alpha_{ij}=0$ \cr \cr Q_{(kl)(ij)}-{|c_{ij}^\pi| \over x_{kl}^0}\ &if $\ c_{ij}^\pi > 0, \ (i,j)\in F, \ \beta_{ij}=0$ \cr\cr  Q_{(kl)(ij)} \ & \ otherwise \cr}$$

\par{\bf Note.} The problem (8) has $2(mn)^2+3(mn)+m+n$ variables which makes it a large problem and hard to solve. A simpler form of IQTP is when the matrix $Q$ is diagonal. In this case the variables would be $5(mn)+m+n$ and it will be easier to handle. Furthermore, the optimal values can be obtain similarly.

\vfill\eject

\par \noindent {\twelvebf 4. Inverse Quadratic Transportation Problem Under ${\rm L_{\infty}}$ norm.}
\bigskip
\noindent Now we want to study the IQTP under ${ L_{\infty}}$ norm. It can be resulted directly form previous section that we just need to change the objective function of (6). Thus, IQTP under $L_{\infty}$ norm is stated as follows:

$${\rm Minimize}\ \| Q-H \|_{\infty} + \| c-d \|_{\infty} \eqno(9a)$$
$$s.t.\ w^1_i+w^1_{n+j}-\sum_{l=1}^m\sum_{k=1}^n H_{(kl)(ij)}x^0_{kl}=d_{ij},\ (i,j)\in F \eqno(9b)$$
$$\ w^1_i+w^1_{n+j}-\sum_{l=1}^m\sum_{k=1}^n H_{(kl)(ij)}x^0_{kl}+w^2_{ij}=d_{ij},\ (i,j)\in L \eqno(9c)$$
$$H_{(kl)(ij)}=H_{(ij)(kl)}, \ (i,j),(k,l)\in J \eqno(9d)$$
$$w^2_{ij}\ge 0,\ (i,j)\in J \eqno(9e)$$

\par We know that $\| c-d \|_{\infty}={\rm max}_{(i,j)\in J} |c_{ij}-d_{ij}|$ and $\| H-Q \|_{\infty}={\rm max}_{(k,l)\in J} (\sum_{(i,j)\in J} $ $|H_{(kl)(ij)}-Q_{(kl)(ij)}|)$. We restate the problem (9) using an identical substitution for the absolute terms and defining $\theta_1={\rm max}_{(k,l)\in J} (\sum_{(i,j)\in J} [\Gamma_{(kl)(ij)}+\Delta_{(kl)(ij)}])$ and $\theta_2={\rm max}_{(i,j)\in J}(\alpha_{ij}+\beta_{ij})$. Hence, the IQTP under $L_{\infty}$ norm is as follows:

$${\rm Minimize}\ \theta_1+ \theta_2 \eqno(10a)$$
$$s.t.\ \sum_{l=1}^m\sum_{k=1}^n [\Delta_{(kl)(ij)}-\Gamma_{(kl)(ij)}]x^0_{kl}-\alpha_{ij}+\beta_{ij}=c^{\pi}_{ij},\ (i,j)\in F \eqno(10b)$$
$$\ \sum_{l=1}^m\sum_{k=1}^n [\Delta_{(kl)(ij)}-\Gamma_{(kl)(ij)}]x^0_{kl}-\alpha_{ij}+\beta_{ij}+w^2_{ij}=c^{\pi}_{ij},\ (i,j)\in L \eqno(10c)$$
$$\Gamma_{(kl)(ij)}=\Gamma_{(ij)(kl)}, \ (i,j),(k,l)\in J \eqno(10d)$$
$$\Delta_{(kl)(ij)}=\Delta_{(ij)(kl)}, \ (i,j),(k,l)\in J \eqno(10e)$$
$$\theta_1 \ge \sum_{j=1}^m\sum_{i=1}^n [\Delta_{(kl)(ij)}+\Gamma_{(kl)(ij)}],\ (i,j)\in F\cup L \eqno(10f)$$
$$\theta_2\ge \alpha_{ij}+\beta_{ij},\ (i,j)\in F\cup L \eqno(10g)$$
$$\Gamma_{(kl)(ij)},\ \Delta_{(kl)(ij)},\ \alpha_{ij},\ \beta_{ij},\ w^2_{ij}\ge 0,\ (i,j),(k,l)\in J \eqno(10h)$$

\noindent where $c^{\pi}_{ij}=c_{ij}-(w^1_i+w^1_{n+j})+\sum_{l=1}^m\sum_{k=1}^n (Q_{(kl)(ij)}x^0_{kl})$.
\par The same as IQTP under $L_1$ norm, by placing some conditions on $c_{ij}^\pi$ and considering that $\alpha_{ij}$ and $\beta_{ij}$ cannot take positive values and similarly $\Delta_{(kl)(ij)}$ and $\Gamma_{(kl)(ij)}$ cannot take positive values, we can obtain the optimal value of IQTP under $L_{\infty}$ norm as follows:
$$d_{ij}^*=\cases {c_{ij}+|c_{ij}^\pi| \ &if $\ c_{ij}^\pi < 0, \ (i,j)\in F\cup L , \ \alpha_{ij}>0$ \cr \cr  c_{ij}-|c_{ij}^\pi| \ &if $\ c_{ij}^\pi > 0, \ (i,j)\in F, \ \beta_{ij}>0$ \cr \cr c_{ij} \ & \ otherwise \cr}$$
\smallskip

$$H_{(kl)(ij)}^*=\cases {Q_{(kl)(ij)}+{|c_{ij}^\pi| \over x_{kl}^0} \ &if $\ c_{ij}^\pi < 0, \ (i,j)\in F\cup L , \ \alpha_{ij}=0$ \cr \cr Q_{(kl)(ij)}-{|c_{ij}^\pi| \over x_{kl}^0}\ &if $\ c_{ij}^\pi > 0, \ (i,j)\in F, \ \beta_{ij}=0$ \cr\cr  Q_{(kl)(ij)} \ & \ otherwise \cr}$$

which is the same as in the case of $L_1$ norm.

\bigskip
\noindent {\twelvebf 5. Conclusion}
\bigskip
\noindent  Inverse optimization has many practical applications. The special case of inverse problem has been studied in this paper. In particular, we obtained the inverse form to quadratic transportation problem under two norms and the optimal values of the problems have been simplified.

\bigskip
\noindent {\twelvebf References}
\bigskip
\item{[A1]} Adlakha, V.; Kowalski, K. On the Quadratic Transportation Problem. {\it Open Journal of Optimization}. (2013) 89--94.
\medskip
\item{[A2]} Ahuja, R. K.; Orlin, J. B. Inverse Optimization. {\it Operations Research} Vol. 49, No. 5, (2001) 771--783.
\medskip
\item{[B1]} Bazaraa, M. S.; Jarvis J. J.; Sherali H. D. {\it Linear Programming and Network Flows}. A John Wiley and Sons, Inc. 4th ed. (2010).
\medskip
\item{[B2]} Bomze, I. M.; Danninger, G. A Finite Algorithm for Solving General Quadratic Problems. {\it Journal of Global Optimization} 4: (1994) 1--16.
\medskip
\item{[B3]} Bretthauer, K. M.; Shetty, B.; Syam, S. A branch and bound algorithm for integer quadratic knapsack problem. {\it ORSA J. Computing}.7, (1995) 109--116.
\medskip
\item{[B4]} Bretthauer, K. M.; Shetty, B.; Syam, S. A projection method for the integer quadratic knapsack problem. {\it J. Opl. Res. Soc.} 47, (1995) 457--462.
\medskip
\item{[C1]} Chaovalitowangse, W.; Pardalas, P. M.; Prokopyev, O. A. A new linearization technique for quadratic 0-1 programming problems. {\it Operation Research Letters}. 32, (2004) 517--522.
\medskip
\item{[C2]} Cosares, S.; Hochbaum, D. S. Strongly Polynomial Algorithms for a Quadratic Transportation Problem with a Fixed Number of Sources. {\it Mathematics of Operation Research} Vol. 19 No. 1 (1994).
\medskip
\item{[J1]} Jain, S.; Arya, N. Inverse Optimization for Quadratic Programming Problems. $\ \ $ {\it IJORN} 2 (2013) 49--56.
\medskip
\item{[K1]} Kozlov, M. K.; Tarasov, S. P.; Khachiyan, L. G. Polynomial solvability of convex quadratic programming. {\it Doklady Akademii Nauk SSSR} 248: (1979) 1049--1051.
\medskip
\item{[Y1]} Yang, C.; Zhang, J.; Ma, Z. Inverse maximum flow and minimum cut problem. {\it Optimization} 40 (1997) 147--170.
\medskip
\item{[X1]} Xu, S.; Zhang, J. An inverse problem of the weighted shortest path problem. {\it Japanese J. Indust. Appl. Math.} 12 (1995) 47--59.
\medskip
\item{[Z1]} Zhang, J.;  Cai, M. Inverse problem of minimum cuts. {\it ZOR Math. Methods Oper. Res.} 48 (1998) 51--58.
\medskip
\item{[Z2]} Zhang, J.; Liu, Z. Calculating some inverse linear programming problems. {\it Journal of Computational and Applied Mathematics}. 72, (1996) 261-273.
\medskip
\item{[Z3]} Zhang, J., Zhang, Li. And Xiao, X., 2010, A perturbation approach for an inverse quadratic programming problem. {\it Mathematical Methods of Operations Research.} Vol. 72, 3, (2010) 379--404

\end